\documentclass{amsart}
\numberwithin{equation}{section}

\usepackage{amsmath,amssymb}

\def\p{\psi}

\def\supp{\text{supp}}

\newtheorem{thm}{Theorem}[section]
\newtheorem{lem}[thm]{Lemma}

\newtheorem{rem}[thm]{Remark}

\def\re{\mathbb{R}}

\begin{document}

\title[An inverse problem of determining the orders ... ]{An inverse problem of determining the orders of systems of fractional pseudo-differential equations}

\author{Ravshan Ashurov}
\author{Sabir Umarov}
\date{}

\address{
$^1$Institute of Mathematics of Academy of Sciences of Republic of Uzbekistan}
\email{ashurovr@gmail.com}

\address{
$^2$University of New Haven, Department of Mathematics, 300 Boston Post Road, \\ West Haven, CT 06516, USA}
 \email{sumarov@newhaven.edu}

\begin{abstract} As it is known various dynamical processes can be modeled through the systems of time-fractional order pseudo-differential equations. In the modeling process one frequently faces with determining adequate orders of time-fractional derivatives in the sense of Riemann-Liouville or Caputo. This problem is qualified as an inverse problem.
The right (vector) order can be found utilizing the available data. In this paper we consider this inverse problem for
linear systems of fractional order pseudo-differential equations. We prove that the Fourier transform of the
vector-solution $\widehat{U}(t, \xi)$ evaluated at a fixed time instance, which becomes possible due to the available data,
recovers uniquely the unknown vector-order of the system of governing pseudo-differential equations.

\end{abstract}

\maketitle

\vspace{2pc}
\noindent
{\it Keywords}: system of differential
equations, fractional order differential equation,
pseudo-differential operator, matrix symbol, inverse problem,
determination of the fractional derivative's order

\section{Introduction}\label{sec:1}

In modern science and engineering researchers frequently use fractional order differential equations for modeling of dynamics of various complex stochastic processes arising in different fields; see, for
example, \cite{Hilfer,Mainardi,MetzlerKlafter00,West} in physics,
\cite{MachadoLopes,SGM} in finance,   \cite{BMR} in hydrology, \cite{Magin} in cell
biology,  among others. In the last few decades  several books,  devoted to fractional order differential and
pseudo-differential equations and their various applications, have been published (see e.g.
\cite{SKM,Podlubny,KST,Handbook,Umarov,UHK_book_2018}).

In fractional order modeling, in contrast to integer order equations, orders of fractional order governing equations are often unknown, and requires to utilize available data to measure. Therefore, one of the key questions arising in the process of modeling is a proper determination of  a fractional order of the governing equation. The problem of determining a correct order of an equation (or orders of equations if the model uses more than one governing equation) is classified as an inverse problem. Inverse problems naturally require additional conditions (or information) for a solution. For subdiffusion equations, in which the order is between zero and one, the inverse problem of determination of the order has been considered by a number
of authors; see  \cite{AA,AF,AU,Che,Jan,LiLiu,LiL,LiY}  and references therein. For the survey paper, we refer the reader to \cite{LiLiu} by Li et. al.  Note that in all the refereed works the subdiffusion equation was considered in a bounded
domain $\Omega\subset \mathbb{R}^N$. In addition, it should be
noted that in publications \cite{Che,Jan,LiL,LiY} the following
relation was taken as an additional condition
\begin{equation}\label{ex1}
u(x_0,t)= h(t), \,\, 0<t<T,
\end{equation}
at a monitoring point $x_0\in \overline{\Omega}$. But this
condition, as a rule (an exception is the work \cite{Jan} by J.
Janno, where both the uniqueness and existence are proved), can
ensure only the uniqueness of the solution of the inverse problem
\cite{Che,LiL,LiY}.
In paper \cite{AU} authors obtained the existence and uniqueness result, considering as an additional condition, the value of the projection of the
solution onto the first eigenfunction of the elliptic part of the
subdiffusion equation. Note that the
technique used in \cite{AU} is applicable only when the first
eigenvalue is zero. In more general case the uniqueness and existence of a solution of the inverse problem of determination of an unknown
order of the fractional derivative in the subdiffusion equation
was proved in the recent work \cite{AA}. In
this case, the additional condition is $ || u (x, t_0) ||^2 = d_0
$, and the boundary condition is not necessarily homogeneous. In paper \cite{AF} authors studied the inverse problem
for the simultaneous determination of the order of the
Riemann-Liouville time fractional derivative and the source
function in the subdiffusion equations.

The purpose of this work is to investigate the inverse problem of
determining the vector-order
of the time-fractional derivatives of systems of pseudo-differential equations. We note that
systems (linear and non-linear) of fractional
ordinary equations and partial differential equations have rich applications and are used in modeling of various processes arising in modern science and engineering.  For example, they are used in modeling of
processes in biosystems \cite{DasGupta,Rihan,GuoFang}, ecology
\cite{Khan,Rana}, epidemiology \cite{Zeb,Islam}, etc.

In this paper we consider the following system of linear
homogenous time-fractional order pseudo-differential equations

    \begin{equation} \label{system_0}
     \begin{cases}
     \mathcal{D}^{\beta_1}u_1(t,x) = A_{1, 1}(D) u_1(t,x) +\dots A_{1, m}(D) u_m(t,x), \\
         \mathcal{D}^{\beta_2}u_2 (t,x)= A_{2, 1}(D) u_1 (t,x)+\dots A_{2, m} (D)u_m(t,x), \\
          {\cdots } \\
         \mathcal{D}^{\beta_m}u_m (t,x) = A_{m, 1}(D) u_1(t,x) +\dots A_{m, m}(D) u_m (t,x).
           \end{cases}
     \end{equation}
where $\mathcal{B}=\langle \beta_1, \dots, \beta_m \rangle,$  $0<\beta_j \le 1, \ j=1,\dots, m,$
is an unknown vector-order to be determined, the operator
$\mathcal{D}$ on the left hand side expresses either the Riemann-Liouville
derivative $D_+$ or the Caputo derivative $D_{\ast}$, and
$A_{j,k}(D)$ are pseudo-differential operators with (possibly
singular) symbols depending only on dual variables (for
simplicity) and described later. The initial conditions depend on
the form of fractional derivatives.

 As it follows from
our main result,
a predetermined value of the Fourier transform $\hat{U}(t,\xi) = \langle \hat{u}_1(t,\xi), \dots, \hat{u}_m(t,\xi) \rangle$ of the
solution $U(t, x)= \langle u_1(t,x), \dots, u_m(t,x) \rangle$ of the initial value problem for system \eqref{system_0} at
an appropriate fixed point $\xi_0\in \mathbb{R}^m$ satisfying some condition (see Eq. \eqref{xi0} ), that is
\begin{equation}\label{3}
\hat{u}_j (t_0, \xi^0)=d_j,\quad j=1,2, ..., m,
\end{equation}
where $t_0 \geq 1$ is an observation time, uniquely recovers the
vector-order $\mathcal{B}=\langle \beta_1, \dots, \beta_m \rangle$
of the fractional derivatives.

 In the particular case, the determining of a scalar order for one equation was considered in \cite{AU} and the forward problem for systems of pseudo-differential equations in \cite{UAY}. From this point of view the current paper is a logical continuation of these two papers.

We note that in \cite{AU} the additional condition is represented
in the form
 \[
 \int\limits_{\Omega} u(t_0, x)v_1(x)dx=d \neq 0,
 \]
that is, in the form of a projection of the solution $u(t_0, x)$
onto the first eigenfunction $v_1(x)$ of the elliptic part of the
equation considered in an arbitrary bounded domain $\Omega\subset
\re^n$. Condition \eqref{3} can be considered as the projection of
the solution $u(t_0, x)$ onto "the eigenfunction" $e^{-ix\xi_0}$:
 \[
 \int\limits_{\re^n} u(t_0, x)e^{-ix\xi_0} dx=d .
 \]

To our best knowledge, the inverse problem for the system of
equations \eqref{system_0} with the vector-order fractional
derivatives under the additional condition \eqref{3} is considered
for the first time.

\section{Main results}

\subsection{Notations.} We follow the notions and notations introduced in \cite{UAY}. For the reader's convenience, below we introduce the main notations used in the current paper; for details see \cite{Umarov,UAY}. Let $G \subseteq \mathbb{R}^n$  be an open set and $p \ge 1.$ The space ${\mathbf \Psi}_{G,p}(\mathbb{R}^n)$ comprises of
functions $\psi\in  L_p(\mathbb{R}^n)$, such that $\supp \
\hat{\psi}\Subset G$, i.e. the Fourier transform
\[
\hat{\psi}(\xi) = \int\limits_{\re^n} f(x) e^{-ix\xi} dx
\]
of $\psi$ has a compact support in $G$.
This is a topological-vector space with respect to the following convergence: a sequence $\p_n \to \p$ if $\supp \, \hat{\p}_n \Subset G,$ and $\p_n \to \p$ in $L_p(\re^n).$ For relations of the spaces ${\mathbf
\Psi}_{G,p}(\mathbb{R}^n)$ to Sobolev spaces and Schwartz
distributions see \cite{Umarov}.

Let $A(\xi)$ be a continuous function in $G$. Outside of $G$ or on
its boundary $A(\xi)$ may have singularities of arbitrary type.
For a function $\varphi \in {\mathbf \Psi}_{G,p}(\mathbb{R}^n)$
the pseudo-differential operator $A(D)$ corresponding to the
symbol $A(\xi)$ is defined  by the formula

\begin{equation}
\label{04} A(D)\varphi (x)= \frac{1}{(2 \pi)^n} \int_G A(\xi)
\hat{\varphi} (\xi) e^{i x \xi} d \xi\, \quad x \in \mathbb{R}^n.
\end{equation}
For the systematic presentation of the theory of
pseudo-differential operators being considered in this paper we
refer the reader to \cite{Umarov}.

Let
\[
\mathbb{A}(D)= \begin{bmatrix}
        A_{1,1}(D) & \dots & A_{1, m} (D) \\
        \dots           & \dots & \dots  \\
        A_{m, 1}(D) & \dots & A_{m,m} (D)
        \end{bmatrix}\]
be the matrix pseudo-differential operator with constant (that is
not depending on the variable $x$) matrix-symbol
\begin{equation} \label{matrix_0}
\mathcal{A}(\xi) = \begin{bmatrix}
        A_{1,1}(\xi) & \dots & A_{1, m} (\xi) \\
        \dots           & \dots & \dots  \\
        A_{m, 1}(\xi) & \dots & A_{m,m} (\xi)
        \end{bmatrix},
\end{equation}
defined and continuous in $G$ in the sense of the matrix norm.

With the matrix form of the pseudo-differential operator we can
represent system (\ref{system_0}) in the vector form:
\begin{equation}
\label{system_1} \mathcal{D}^{\mathcal{B}} {U}(t,x) =
\mathbb{A}(D) {U} (t,x),
\end{equation}
where $ \mathcal{D}^{\mathcal{B}} {U}(t,x)=\langle
\mathcal{D}^{{\beta_1}} {u}_1(t,x), \dots, \mathcal{D}^{{\beta_m}}
{u}_m(t,x) \rangle. $

Below we will use two main forms of fractional derivatives, namely the Riemann-Liouville form and the Caputo form. Let $k$ be a natural number and $k-1 \le \beta <k.$ Then the
fractional derivative of order $\beta$ of a measurable function
$f$ in the sense of Riemann--Liouville is defined as
\[
D^{\beta}_{+} f(t) = \frac{1}{\Gamma(k-\beta)}\frac{d^k}{dt^k}
\int_0^t \frac{f(\tau)d\tau}{(t-\tau)^{\beta+1-k}},
\]
provided the expression on the right exists. Here $\Gamma(t)$ is
Euler's gamma function. If we replace differentiation and
fractional integration in this definition, then we get the
definition of a regularized derivative, that is, the definition of
a fractional derivative in the sense of Caputo:
\[
D^{\beta}_{\ast} f(t) = \frac{1}{\Gamma(k-\beta)} \int_0^t
\frac{f^{(k)}(\tau)d\tau}{(t-\tau)^{\beta+1-k}},
\]
provided the integral on the right exists.

We assume that the matrix-symbol is symmetric,
$A_{k,j}(\xi)=A_{j,k}(\xi)$ for all $k,j=1,\dots,m,$ and $\xi\in
G,$ and diagonalizable. Namely, there exists an invertible
$(m\times m)$-matrix-function $M(\xi),$ such that
\begin{equation}
\label{matrix_01} \mathcal{A}(\xi) = M^{-1}(\xi) \Lambda(\xi)
M(\xi), \quad \xi \in G,
\end{equation}
with a diagonal matrix

\begin{equation} \label{matrix_02}
{\Lambda}(\xi) = \begin{bmatrix}
        \lambda_1(\xi) & \dots & 0  \\
        \dots           & \dots & \dots  \\
        0  & \dots & \lambda_m (\xi)
        \end{bmatrix}.
\end{equation}
We denote entries of matrices $M(\xi)$ and $M^{-1}(\xi)$ by
$\mu_{j,k}(\xi), \ j,k=1,\dots,m,$ and $\nu_{j,k}(\xi)$,
$j,k=1,\dots,m,$ respectively.

Since initial conditions depend on the form of the fractional
derivative on the left hand side of equation (\ref{system_1}), we
will consider the cases with the Caputo and Riemann-Liouville
derivatives separately. We first formulate our main result in the
case of Caputo fractional derivative. The case of Rieman-Liouville
fractional derivative can be treated similarly.

\subsection{Forward problem}
Let $\mathcal{B}$ be a known vector-order with $0<\beta_j\leq 1$,
$j=1, \dots, m$. Consider the following Cauchy problem

\begin{align} \label{Cauchy_01_h}
D_{\ast}^{\mathcal{B}} {U}(t,x) &= \mathbb{A}(D) {U} (t,x), \quad
t>0, \ x \in \mathbb{R}^n,
\\
U(0,x) & =\varPhi(x), \quad x \in \mathbb{R}^n,
\label{Cauchy_02_h}
\end{align}
where $\varPhi(x) = \langle \varphi_1(x), \dots, \varphi_m(x)
\rangle\in {\mathbf \Psi}_{G,p}(\mathbb{R}^n)$ and the fractional
derivatives on the left are in the sense of Caputo.

We call the Cauchy problem (\ref{Cauchy_01_h})-(\ref{Cauchy_02_h})
\emph{the forward problem.}

A representation formula for the solution of the forward problem
was obtained in \cite{UAY} and it has the form
\[
u_j(t,x)= \frac{1}{(2\pi)^n} \int\limits_{\mathbb{R}^n}\sum\limits_{k=1}^m
s_{j,k}(t,\xi)\hat{\varphi}_k(\xi)e^{ix\xi}d\xi, \quad j=1, \dots,
m,
\]
where
\[
s_{j,k}(t,\xi)=\sum\limits_{l=1}^m
\mu_{j,l}(\xi)\nu_{l,k}(\xi)E_{\beta_l}(\lambda_l(\xi)
t^{\beta_l}).
\]
Here we denoted by $E_{\beta_j}(z), j=1,\dots,m,$ the
Mittag-Leffler functions of indices $\beta_1,\dots, \beta_m,$ respectively.

We rewrite function $u_j$ as
\begin{align*}
u_j(t,x) &= \frac{1}{(2\pi)^n} \int\limits_{\mathbb{R}^n}\sum\limits_{k=1}^m
\sum\limits_{l=1}^m
\mu_{j,l}(\xi)\nu_{l,k}(\xi)E_{\beta_l}(\lambda_l(\xi)
t^{\beta_l})\hat{\varphi}_k(\xi)e^{ix\xi}d\xi
\\
&= \frac{1}{(2\pi)^n} \int\limits_{\mathbb{R}^n}\sum\limits_{l=1}^m
E_{\beta_l}(\lambda_l(\xi)
t^{\beta_l})\bigg[\mu_{j,l}(\xi)\sum\limits_{k=1}^m\nu_{l,k}(\xi)\hat{\varphi}_k(\xi)\bigg]e^{ix\xi}d\xi
\\
&= \frac{1}{(2\pi)^n} \int\limits_{\mathbb{R}^n}\sum\limits_{l=1}^m
E_{\beta_l}(\lambda_l(\xi) t^{\beta_l})K_{j,l}\big(\xi,
\widehat{\Phi}(\xi)\big)e^{ix\xi}d\xi,
\end{align*}
where
\[
K_{j,l}\big(\xi,
\widehat{\Phi}(\xi)\big)=\mu_{j,l}(\xi)\sum\limits_{k=1}^m\nu_{l,k}(\xi)\hat{\varphi}_k(\xi),
\]
and
\[
\widehat{\Phi}(\xi)=\langle\hat{\varphi}_1(\xi),\hat{\varphi}_2(\xi),\cdot\cdot\cdot,\hat{\varphi}_m(\xi)\rangle.
\]

For the Fourier transform of the solution we have
\begin{equation}\label{1}
\hat{u}_j(t,\xi)=\sum\limits_{l=1}^m E_{\beta_l}(\lambda_l(\xi)
t^{\beta_l})K_{j,l}\big(\xi, \widehat{\Phi}(\xi)\big).
\end{equation}

Note that under the above conditions on the matrix-symbol
$A_{j,k}(\xi)$ and on the function $\Phi(x)$, this Fourier
transform exists at each point $\xi\in \mathbb{R}^n$.

\subsection{Inverse problem.}
Now let the parameter $\mathcal{B}$ be an unknown vector-order of
the time derivative with $\beta_0\leq\beta_j< 1$, $j=1, \dots, m$,
$\beta_0\in (0,1)$. The main purpose of this paper is to
investigate the inverse problem of identifying of these parameters
$\beta_j$. Since there are $m$ unknown parameters, we need to set
$m$ conditions. We pass on to the determining of these additional
conditions.

In what follows, we will assume that
\begin{equation}\label{lambda}
|\arg \lambda_j(\xi)|> \frac{\pi}{2}, \quad \xi \in G,\quad j=1,
\dots, m.
\end{equation}
Let $\xi^0=(\xi^0_1, \xi^0_2, ..., \xi^0_m)\in G\subset
\mathbb{R}^n$ be a vector such that the determinant of the
matrix
\begin{equation}
\label{matrixK}
\mathcal{K}(\xi^0) \equiv \{K_{j,l}\big(\xi^0, \widehat{\Phi}(\xi^0)\big)\}, \ j, l =1,\dots,m,
\end{equation}
 satisfies the condition
\begin{equation}\label{xi0}
\big|K_{j,l}\big(\xi^0, \widehat{\Phi}(\xi^0)\big)\big|\neq 0.
\end{equation}

To find the unknown parameters $\beta_l, l=1,\dots, m$, we
consider the following additional conditions
\begin{equation}\label{ad1}
f_j(\mathcal{B}, t_0, \xi^0)\equiv \hat{u}_j (t_0,
\xi^0)=d_j,\quad j=1,\dots, m,
\end{equation}
where $d_j$ are given numbers and $t_0$ is defined later (see Lemmas \ref{ECaputo} and \ref{ERL}). We call
problem (\ref{Cauchy_01_h})--(\ref{Cauchy_02_h}) together with the
additional condition (\ref{ad1}) \emph{the inverse problem.}

It follows from (\ref{1}) and (\ref{ad1}) that for all $j=1,2, ...  , m$
\begin{equation}\label{4}
\sum\limits_{l=1}^m E_{\beta_l}(\lambda_l(\xi^0)
t_0^{\beta_l})K_{j,l}\big(\xi^0, \widehat{\Phi}(\xi^0)\big)=d_j.
\end{equation}
These are in fact the system of equations to define the orders
$\beta_l, l=1,2,..., m$.

Due to condition (\ref{xi0}) one can solve system (\ref{4}) with
respect to the Mittag-Leffler functions $E_{\beta_l}$, i.e.
\begin{equation}\label{ad2}
E_{\beta_l}(\lambda_l(\xi^0) t_0^{\beta_l})= b_l, \quad l=1,2,...,
m,
\end{equation}
where $b_l, \, l=1,2,..., m,$ are components of the vector
$\mathcal{K}^{-1}(\xi^0)
\textbf{d},$ with $\textbf{d}=\langle d_1,\dots,d_m \rangle$ and
$\mathcal{K}^{-1}(\xi^0)$ is the inverse matrix to
$\mathcal{K}(\xi^0),$ defined in \eqref{matrixK}. Thus to define
each unknown parameter $\beta_l$ we obtained a separate equation
(\ref{ad2}).

Let $R_{C, l}$ be the range of values of the function $e_{1,
\lambda_l} (\beta_l)\equiv E_{\beta_l}(\lambda_l(\xi^0)
t_0^{\beta_l})$ when $\beta_l$ runs over the half-interval
$[\beta_0, 1)$, i.e.
\[
e_{1, \lambda_l}: [\beta_0, 1)\rightarrow R_{C, l}\subset C,
\]
where $C$ is a complex plane, and the index $C$ emphasizes that we
are considering the case of the Caputo derivatives.

Obviously, for equations (\ref{ad2}), in order to have solutions
with respect to $\beta_l,$ the right-hand sides of these equations
must lie within the values of the functions on the left-hand sides
of these equations, i.e.
\begin{equation}\label{Crange}
b_l\in R_{C, l}, \quad l=1,2,..., m.
\end{equation}

On the other hand, by virtue of Rolle's theorem, the strict
monotonicity  of either the real part or the imaginary
part of the function $E_{\beta_l}(\lambda_l(\xi^0) t_0^{\beta_l})$ in the variable $\beta_l$
is sufficient for the uniqueness of the solution to equation
(\ref{ad2}).

Let us introduce the following notation
\[
R_C(\beta_l)=\Re(E_{\beta_l}(\lambda_l(\xi^0) t_0^{\beta_l})),
\]
where $\Re(z)$ is the real part of $z$ and the index $C$ again
emphasizes that we are considering the case of the Caputo
derivatives. The necessity of condition (\ref{lambda}) is that its fulfillment
guarantees the strict monotonicity of the function $R_C(\beta_l)$
it the variable $\beta_l$ for each fixed $l.$

\subsection{Main results}
\par
The main results of this paper are stated in Theorems \ref{thm_Caputo} and \ref{thm_Riemann-Liouville}.

\begin{thm}
\label{thm_Caputo} Let $\xi^0$ satisfy condition \eqref{xi0} and
$t_0 > T_0,$ where $T_0$ is identified in Lemma \ref{ECaputo}. Let
the numbers $d_l$ on the right hand side of equation (\ref{ad1})
be such that the components $b_l, \ l=1,\dots, m,$ of the vector
$\mathcal{K}^{-1}(\xi^0) \textbf{d}$ satisfy the conditions
(\ref{Crange}). Then for each $l$ there exists the unique number
$\beta_l^\ast\in [\beta_0, 1]$ such that the Fourier transform of
the solution $u_j(t, x), j=1,\dots,m,$  of the forward problem
with $\beta_j=\beta_j^\ast, j=1,\dots,m,$ satisfies equation
(\ref{ad1}).
\end{thm}

The proof of this theorem follows from the existence and uniqueness theorem for the forward problem proved in \cite{UAY} (see Theorem 3.1) and Lemma \ref{ECaputo} below. Therefore, in order to prove the theorem we need to prove only this lemma.
The proof of Lemma \ref{ECaputo} is given in Section \ref{sec_3}.

\begin{lem}
\label{ECaputo}

Given $\beta_0$ in the interval $0<\beta_0< 1$ and $\xi_0$
satisfying condition (\ref{xi0}), there exists a number
$T_0=T_0(\xi^0, \beta_0)>1$, such that for all $t_0\geq T_0$ the
function $R_C(\beta_l)$ is positive and strictly monotonically
decreasing with respect to $\beta_l\in [\beta_0, 1]$ and
\begin{equation}\label{RICE}
R_C(1) \leq R_C(\beta_l) \leq R_C(\beta_0),\,\, l=1, \dots, m.
\end{equation}

\end{lem}

\begin{rem}\label{Uniqueness}
Theorem \ref{thm_Caputo} defines the vector-order $\mathcal{B}^\ast=(\beta_1^\ast,
\beta_2^\ast,..., \beta_m^\ast)$ uniquely from conditions (\ref{ad1}).
Hence, if we define $f_j(\mathcal{B}, \cdot, \cdot)$ at another
time instant $t_1$ and point $\xi^1$ and get a new
$\mathcal{B^{\ast \ast}}$, i.e. $f_j(\mathcal{B}^{\ast \ast}, t_1,
\xi^1) = d_j^1$, then from the equality $f_j(\mathcal{B}^{\ast
\ast}, t_0, \xi^0) = d_j$, by virtue of the theorem, we obtain
$\mathcal{B}^{\ast \ast}=\mathcal{B}^{\ast}$.
\end{rem}

Now consider the following initial-value problem
\begin{align} \label{Cauchy_10}
D_{+}^{\mathcal{B}} {U}(t,x) &= \mathbb{A}(D) {U} (t,x), \quad
t>0, \ x \in \mathbb{R}^n,
\\
J^{1-\mathcal{B}}U(0,x) & =\varPhi(x), \quad x \in \mathbb{R}^n,
\label{Cauchy_20}
\end{align}
where $\varPhi(x) = \langle \varphi_1(x), \dots, \varphi_m(x)
\rangle\in {\mathbf \Psi}_{G,p}(\mathbb{R}^n)$ and  the fractional
derivatives on the left hand side of equation \eqref{Cauchy_10}
are in the sense of Riemann-Liouville.

We call the Cauchy problem (\ref{Cauchy_10})-(\ref{Cauchy_20})
\emph{the second forward problem.}

A representation formula for the solution of the second forward
problem was  also obtained in \cite{UAY} and it has the form
\[
u_j(t,x)= \frac{1}{(2\pi)^n} \int\limits_{\mathbb{R}^n}\sum\limits_{k=1}^m
s^+_{j,k}(t,\xi)\hat{\varphi}_k(\xi)e^{ix\xi}d\xi, \quad j=1,
\dots, m,
\]
where
\[
s^+_{j,k}(t,\xi)=\sum\limits_{l=1}^m
\mu_{j,l}(\xi)\nu_{l,k}(\xi)t^{\beta_l-1}E_{\beta_l,
\beta_l}(\lambda_l(\xi) t^{\beta_l}).
\]
Here we denoted by $E_{\beta_j, \beta_j}(z), j=1,\dots,m,$ the
two-parametric Mittag-Leffler functions.

For the Fourier transform of the solution we have
\begin{equation}\label{Fu+}
\hat{u}_j(t,\xi)=\sum\limits_{l=1}^m t^{\beta_l-1}E_{\beta_l,
\beta_l}(\lambda_l(\xi) t^{\beta_l})K_{j,l}\big(\xi,
\widehat{\Phi}(\xi)\big).
\end{equation}

Suppose that condition (\ref{lambda}) is fulfilled and choose
$\xi^0\in G$ so that inequality (\ref{xi0}) holds.

We call problem (\ref{Cauchy_10})-(\ref{Cauchy_20}) together with
the additional condition (\ref{ad1}) \emph{the second inverse
problem}.

Note that additional condition (\ref{ad1}) is in fact the equation to
determine the unknown parameters $\beta_l$. Performing similar
calculations as above, by virtue of condition (\ref{xi0}), we
rewrite (\ref{ad1}) as
\begin{equation}\label{ad3}
t^{\beta_l-1}E_{\beta_l, \beta_l}(\lambda_l(\xi) t^{\beta_l})=
b_l, \quad l=1,2,..., m,
\end{equation}
where $b_l, \, l=1,2,..., m,$ are the same numbers as above.

Let $R_{RL, l}$ be the range of values of the left-hand side of
these equations when $\beta_l$ runs over the half-interval
$[\beta_0, 1)$.  Here the index $RL$ emphasizes that we are
considering the case of the Riemann-Liouville derivatives.

Again, as in case of equations (\ref{ad2}), a necessary condition
for the existence of solutions to equations (\ref{ad3}) is the
inclusion
\begin{equation}\label{RLrange}
b_l\in R_{RL, l}, \quad l=1,2,..., m.
\end{equation}

On the other hand, by virtue of Rolle's theorem, the strict
monotonicity  of the function
$\Re(t^{\beta_l-1}E_{\beta_l, \beta_l}(\lambda_l(\xi)
t^{\beta_l}))$ in the variable $\beta_l$ is sufficient for the uniqueness of the solution to
equation (\ref{ad3}).

However, if $|\Re(\lambda_l(\xi^0))| = |\Im(\lambda_l(\xi^0))|$
(note, under the condition (\ref{lambda}) one has
$\Re(\lambda_l(\xi^0))<0$), then the principal part of
$\Re(t^{\beta_l-1}E_{\beta_l, \beta_l}(\lambda_l(\xi)
t^{\beta_l}))$ vanishes, and in this case it is necessary to go to
its next term in the asymptotic. Therefore, to simplify the
presentation, we further assume that
\begin{equation}\label{lambda1}
|\Re(\lambda_l(\xi^0))| \neq |\Im(\lambda_l(\xi^0))|
\end{equation}

Let us introduce the following notation
$$
R_{RL}(\beta_l)=sign(|\Re(\lambda_l(\xi^0))| -
|\Im(\lambda_l(\xi^0))|)\Re(t^{\beta_l-1}E_{\beta_l,
\beta_l}(\lambda_l(\xi) t^{\beta_l})).
$$
Here the index $RL$  again emphasizes that we are considering the
case of the Riemann-Liouville derivatives.

\begin{lem}\label{ERL}
 Given $\beta_0$ in the interval $0<\beta_0< 1$ and $\xi_0$ satisfying condition (\ref{xi0}),
 there exists a number $T_1=T_1(\xi^0,
\beta_0)>1$, such that for all $t_0\geq T_1$ the function
$R_{RL}(\beta_l)$ is positive and strictly monotonically
decreasing with respect to $\beta_l\in [\beta_0, 1]$ and
\begin{equation}\label{RICE_1}
R_{RL}(1) \leq R_{RL}(\beta_l) \leq R_{RL}(\beta_0),\,\, l=1,
\dots, m.
\end{equation}

\end{lem}

This lemma, similar to the Caputo derivative case, immediately
implies the following main result of this paper in the case of the
Riemann-Liouville  derivatives. The existence and uniqueness
theorem of the corresponding forward problem is proved in
\cite{UAY} (see Theorem 3.4). The proof of Lemma \ref{ERL} is
presented in the next section.

\begin{thm}
\label{thm_Riemann-Liouville} Let $\xi^0$ satisfy condition
\eqref{xi0} and $t_0>T_1,$ where $T_1$ is identified in Lemma
\ref{ERL}. Let the numbers $d_l$ from condition (\ref{ad1}) be
such that the corresponding numbers $b_l$ satisfy the conditions
(\ref{RLrange}). Then for each $l$ there exists the unique number
$\beta_l^\ast\in [\beta_0, 1]$ such that the Fourier transform of
the solution $u_j(t, x)$ of the second forward problem with
$\beta_j=\beta_j^\ast$ satisfies the equation (\ref{ad1}).
\end{thm}

Similar to the Caputo derivative case, Theorem \ref{thm_Riemann-Liouville} defines the vector-order $\mathcal{B}^\ast=(\beta_1^\ast,
\beta_2^\ast,..., \beta_m^\ast)$ uniquely from conditions (\ref{ad1}); see Remark \ref{Uniqueness}.


\section{Proofs of Lemmata \ref{ECaputo} and \ref{ERL}}
\label{sec_3}

Let us denote by $\delta(1; \theta)$ a contour oriented by
non-decreasing $\arg \zeta$ consisting of the following parts: the
ray $\arg  \zeta = -\theta$ with $|\zeta|\geq 1$, the arc
$-\theta\leq \arg \zeta \leq \theta$, $|\zeta|=1$, and the ray
$\arg \zeta = \theta$, $|\zeta|\geq 1$. If $0<\theta <\pi$, then
the contour $\delta(1; \theta)$ divides the complex $\zeta$-plane
into two unbounded parts, namely $G^{(-)}(1;\theta)$ to the left
of $\delta(1; \theta)$ by orientation, and $G^{(+)}(1;\theta)$ to
the right of it. The contour $\delta(1; \theta)$ is called the
Hankel path.

In what follows, we fix $l$ out of $1,\dots, m$ and denote
$\lambda=\lambda_l(\xi^0)$ and $\rho=\beta_l$. Let
$\lambda=-\lambda_1+i\lambda_2$ and by virtue of condition
(\ref{lambda}) one has $\lambda_1>0$.
 Further let
$\theta = (\frac{\pi}{2}+\varepsilon)\rho$, $\alpha =
(\frac{\pi}{2}+2\varepsilon)\rho$, $\rho\in [\beta_0, 1)$ and
$\varepsilon>0$ be such that
$\varepsilon\equiv\varepsilon(\xi^0)<\frac{1}{2}\min\{|\arg
\lambda(\xi^0)|-\pi/2, \pi/2\}$.

Then
\[
\frac{\pi}{2}\rho<\theta<\alpha< \pi \rho, \quad \alpha <|\arg
\lambda|,
\]
and therefore $\lambda t_0^\rho\in G^{(-)}(1;\theta)$.

\

\subsection{Proof of Lemma \ref{ECaputo}}
First we prove Lemma \ref{ECaputo}. By the definition of contour $\delta(1; \theta)$, we have
(see \cite{Dzh66}, formula (2.29), p. 135)
\begin{equation}\label{E1rho}
e_{1,\lambda}(\rho)\equiv E_{\rho}(\lambda t_0^\rho)=
-\frac{1}{\lambda t_0^\rho \Gamma (1-\rho)}+\frac{1}{2\pi i
\rho\lambda
t_0^\rho}\int\limits_{\delta(1;\theta)}\frac{e^{\zeta^{1/\rho}}\zeta}{\zeta+\lambda
t_0^\rho} d\zeta = f_1(\rho)+f_2(\rho).
\end{equation}

To prove the lemma, we need to determine the sign of the real part
of the derivative $\frac{d}{d\rho} e_{1,\lambda} (\rho)$. It is
not hard to estimate the derivative $f'_1(\rho)$. Indeed, let
$\Psi(\rho)$ be the logarithmic derivative of the gamma function
$\Gamma(\rho)$ (for the definition and properties of $\Psi$ see
\cite{Bat}). Then $\Gamma'(\rho) = \Gamma (\rho) \Psi(\rho)$, and
therefore,
$$
f_1'(\rho)=\frac{\ln t_0 - \Psi (1-\rho)}{\lambda t_0^\rho \Gamma
(1-\rho)}.
$$
Since
\[
\frac{1}{\Gamma(1-\rho)}=\frac{1-\rho}{\Gamma(2-\rho)}, \quad
\Psi(1-\rho)=\Psi(2-\rho)-\frac{1}{1-\rho},
\]
the function $f_1'(\rho)$  can be represented as follows
\[
f_1'(\rho)=\frac{1}{\lambda t_0^\rho} \frac{(1-\rho)[\ln t_0 -
\Psi (2-\rho)]+1}{ \Gamma (2-\rho)}.
\]

If $\gamma\approx 0,57722$ is the Euler-Mascheroni constant, then
$-\gamma <\Psi(2-\rho)< 1-\gamma$. By virtue of this estimate we
may write
\begin{equation}\label{Rf1}
-\Re(f_1'(\rho))\geq
\frac{\lambda_1}{|\lambda|^2}\frac{(1-\rho)[\ln t_0 -
(1-\gamma)]+1}{\Gamma(2-\rho) t_0^\rho}\geq
\frac{\lambda_1}{|\lambda|^2 t_0^{\rho}},
\end{equation}
provided
\begin{equation}
\label{new_t0} \ln t_0 > 1-\gamma \quad \text{or} \quad t_0> T_0=e^{1-\gamma}>1.
\end{equation}

To estimate the derivative $f'_2(\rho)$, we denote the integrand
in (\ref{E1rho}) by $F(\zeta, \rho)$:
\[
F(\zeta, \rho)=\frac{1}{2\pi i \rho\lambda t_0^\rho}\cdot
\frac{e^{\zeta^{1/\rho}}\zeta}{\zeta+\lambda t_0^\rho}.
\]
Note, that the domain of integration $\delta(1; \theta)$ also
depends on $\rho$. To take this circumstance into account when
differentiating the function $f'_2(\rho)$, we rewrite the integral
(\ref{E1rho}) in the form:
\[
f_2(\rho)=f_{2+}(\rho)+f_{2-}(\rho)+f_{21}(\rho),
\]
where
\[
f_{2\pm}(\rho)=e^{\pm i \theta}\int\limits_1^\infty F(s\,e^{\pm i
\theta}, \rho)\, ds,
\]

\[
f_{21}(\rho) = i \int\limits_{-\theta}^{\theta} F(e^{i y}, \rho)\,
e^{iy} dy= i\theta \int\limits_{-1}^{1} F(e^{i \theta s}, \rho)\,
e^{i\theta s} ds.
\]

Let us consider the function $f_{2+}(\rho)$. Since $\theta =
(\frac{\pi}{2}+\varepsilon)\rho$ and $\zeta= s\, e^{i\theta}$,
then
\[
e^{\zeta^{1/\rho}}=e^{-s^{\frac{1}{\rho}}(\varepsilon_1-i\varepsilon_2)},
\,\, \cos(\frac{\pi}{2}+\varepsilon) =-\varepsilon_1<0.\,\,
\sin(\frac{\pi}{2}+\varepsilon) =\varepsilon_2>0.
\]
The derivative of the function $f_{2+}(\rho)$ has the form
$$
f_{2+}'(\rho)=\frac{1}{2\pi i \rho\lambda
t_0^\rho}\int\limits_1^\infty
\frac{e^{-s^{\frac{1}{\rho}}(\varepsilon_1-i\varepsilon_2)}\,s\,e^{2ia\rho}\mathcal{M}(s)}{s\,e^{ia\rho}+\lambda t_0^\rho} ds,
$$
where $a=\frac{\pi}{2}+\varepsilon,$ and
\[ \mathcal{M}(s)=-\frac{\varepsilon_1-i\varepsilon_2}{\rho^2}s^{1/\rho}\ln
s+2ia -\frac{1}{\rho}-\ln t_0-\frac{ia s\,e^{ia\rho}+\lambda
t_0^\rho \ln t_0}{s\,e^{ia\rho}+\lambda
t_0^\rho}.
\]
It is not hard to verify, that 
\[
|s\,e^{ia\rho}+\lambda
t_0^\rho|\geq |\lambda| t_0^\rho \sin(\alpha - \theta)\geq
\frac{2}{\pi}|\lambda| t_0^\rho \varepsilon.
\]
 Therefor we arrive
at
$$
|f_{2+}'(\rho)|\leq \frac{C}{\rho(\varepsilon|\lambda|
t_0^\rho)^2}\int\limits_1^\infty
e^{-\varepsilon_1\,s^{1/\rho}}s\,\left[\frac{1}{\rho^2}s^{1/\rho}\ln
s+\ln t_0\right] ds,
$$
or
$$
|f_{2+}'(\rho)|\leq \frac{C}{(|\lambda|
t_0^\rho)^2}\left[\frac{1}{\rho}+\ln t_0\right],
$$
where the constant $C$ depends only on $\varepsilon$ (and
therefore only on $\xi^0$).

The function $f'_{2-}(\rho)$ has exactly the same estimate.

Now consider the function $f_{21}(\rho)$. It is not hard to verify
that
\[
f'_{21}(\rho)=\frac{a}{2\pi \lambda
t_0^\rho}\int\limits_{-1}^1\frac{e^{e^{ias}}e^{2ia\rho
s}\big[2ias-\ln t_0-\frac{ias e^{ia\rho s}+\lambda t_0^\rho \ln
t_0}{e^{ia\rho s}+\lambda t_0^\rho}\big]}{e^{ia\rho s}+\lambda
t_0^\rho} ds.
\]
Therefore,
\[
|f'_{21}(\rho)|\leq C\frac{\ln t_0}{(|\lambda| t_0^\rho)^2}.
\]
Now we show that the real part of the derivative $\frac{d}{d\rho}
e_{1,\lambda}(\rho)$ is negative.
Taking into account estimate (\ref{Rf1}) and the estimates for
$f'_{2\pm}$ and $f'_{21}$, we have
\begin{equation}\label{der}
\Re\bigg(\frac{d}{d\rho} e_{1,\lambda}(\rho)\bigg)<
-\frac{\lambda_1}{|\lambda|^2 t_0^{\rho}}+C\frac{1/\rho+\ln
t_0}{(|\lambda| t_0^{\rho})^2}.
\end{equation}
In other words, this derivative is negative if
\[
t_0^{\rho}>C\frac{1/\rho+ \ln t_0}{\lambda_1}
\]
for all $\rho\in [\beta_0, 1).$ Hence
\begin{equation}\label{t01}
t_0^{\beta_0}>C\frac{1/{\beta_0}+ \ln t_0}{\lambda_1}.
\end{equation}

Thus, there exists a number $T_0 = T_0(\xi^0, \beta_0)>1$ (see
\eqref{new_t0}) such, that for all $t_0\geq T_0$ we have the estimate
\[
\Re\bigg(\frac{d}{d\rho} e_{1,\lambda}(\rho)\bigg)<
0\,\,\text{for}\,\,\text{all} \quad \rho\in [\beta_0,1].
\]

The positivity of $R_C(\beta_l)$ follows from the explicit form of
the function $f_1(\rho)$.

Lemma \ref{ECaputo}, and therefore Theorem \ref{thm_Caputo} are
completely proved.

\

\subsection{Proof of Lemma \ref{ERL}} We now turn to the proof of Lemma \ref{ERL}.
By the definition of contour $\delta(1; \theta)$, we have for
$e_{2,\lambda}(\rho)\equiv t_0^{\rho-1}E_{\rho, \rho}(\lambda
t_0^\rho)$ the following equation  (see \cite{Dzh66}, formula
(2.29), p. 135)
\begin{equation}\label{Erho}
e_{2,\lambda}(\rho)= - \frac{1}{\lambda^2 t_0^{\rho+1} \Gamma
(-\rho)}+\frac{1}{2\pi i \rho \lambda^2
t_0^{\rho+1}}\int\limits_{\delta(1;\theta)}\frac{e^{\zeta^{1/\rho}}\zeta^{\frac{1}{\rho}+1}}{\zeta+\lambda
t_0^\rho} d\zeta = g_1(\rho)+g_2(\rho).
\end{equation}

Since the positivity of  $R_{RL}(1)$ is obvious, then in order to
prove Lemma \ref{ERL} it suffices to show that the derivatives of
$R_{RL}(\rho)$ is negative for all $\rho\in [\beta_0, 1)$.

For the derivative $g'_1(\rho)$ we have
$$
g_1'(\rho)=\frac{\ln t_0 - \Psi (-\rho)}{\lambda^2 t_0^{\rho
+1}\Gamma (-\rho)}.
$$
To get rid of the singularity in the denominators, we use the
equalities
\begin{align*}
\frac{1}{\Gamma(-\rho)} &= -\frac{\rho}{\Gamma(1-\rho)}=-\frac{\rho(1-\rho)}{\Gamma(2-\rho)},
\\
\Psi(-\rho) &=\Psi(1-\rho)
+\frac{1}{\rho}=\Psi(2-\rho)+\frac{1}{\rho}-\frac{1}{1-\rho}.
\end{align*}
Then the function $g_1'(\rho)$  can be represented as follows
\begin{equation}\label{f1}
g_1'(\rho)=\frac{1}{\lambda^2 t_0^{\rho+1}}
\frac{\rho(1-\rho)[\Psi (2-\rho)-\ln t_0]+1-2\rho}{ \Gamma
(2-\rho)}=-\frac{g_{11}(\rho)}{\lambda^2 t_0^{\rho+1}\Gamma
(2-\rho)}.
\end{equation}
Since $\Psi(2-\rho)< 1-\gamma$, then
\[
g_{11}(\rho)>\rho (1-\rho)[\ln t_0 -(1-\gamma)])+2\rho-1.
\]
For $t_0=e^{1-\gamma} e^{2/\rho}$ one has $\rho (1-\rho)[\ln t_0
-(1-\gamma)])+2\rho-1=1$. Hence, $g_{11}(\rho)\geq 1$, provided
$t_0\geq T_1,$ where
\begin{equation}\label{t0}
T_1=  e^{1-\gamma} e^{2/\beta_0} > e^{3-\gamma}>1.
\end{equation}
Thus, by virtue of (\ref{f1}), for
all such $t_0$ we arrive at
\begin{equation}\label{f11}
sign(\lambda_1^2-\lambda_2^2)\Re(g_1'(\rho))\leq-\frac{|\lambda_1^2-\lambda_2^2|}{|\lambda|^4
t_0^{\rho+1}}.
\end{equation}

To estimate the derivative $g'_2(\rho)$, we denote the integrand
in (\ref{Erho}) by $G(\zeta, \rho)$:
\[
G(\zeta, \rho)=\frac{1}{2\pi i \rho\lambda^2 t_0^{\rho+1}}\cdot
\frac{e^{\zeta^{1/\rho}}\zeta^{1/\rho+1}}{\zeta+\lambda t_0^\rho},
\]
and rewrite the integral (\ref{Erho}) in the form:
\[
g_2(\rho)=g_{2+}(\rho)+g_{2-}(\rho)+g_{21}(\rho),
\]
where
\[
g_{2\pm}(\rho)=e^{\pm i \theta}\int\limits_1^\infty G(s\,e^{\pm i
\theta}, \rho)\, ds,
\]

\[
g_{21}(\rho) = i \int\limits_{-\theta}^{\theta} F(e^{i y}, \rho)\,
e^{iy} dy= i\theta \int\limits_{-1}^{1} G(e^{i \theta s}, \rho)\,
e^{i\theta s} ds.
\]

The derivative of the function $ g_{2+}(\rho)$ has the form
$$
g_{2+}'(\rho)=I\cdot\int\limits_1^\infty \frac{e^{-s^{1/\rho}\,
(\varepsilon_1-i\varepsilon_2)}s^{\frac{1}{\rho}+1}\,e^{2ia\rho} \mathcal{N}(s) }{s\,e^{ia\rho}+\lambda t_0^\rho} ds,
$$
where $I=e^{ia}(2\pi i \rho\lambda^2 t_0^{\rho+1})^{-1},$
$a=\frac{3\pi}{4},$ and
\[
\mathcal{N}(s) = -\frac{1}{\rho^2}((\varepsilon_1+i\varepsilon_2)s^{\frac{1}{\rho}}+1)\ln
s+2ia -\frac{1}{\rho}-\ln t_0-\frac{ia s\,e^{ia\rho}+\lambda
t_0^\rho \ln t_0}{s\,e^{ia\rho}+\lambda
t_0^\rho}.
\]
By virtue of the inequality
$|s\,e^{ia\rho}+\lambda t_0^\rho|\geq \frac{2}{\pi}|\lambda|
t_0^\rho \varepsilon$ we arrive at
$$
|g_{2+}'(\rho)|\leq \frac{C}{\rho |\lambda|^3
t_0^{2\rho+1}}\int\limits_1^\infty
e^{-\varepsilon_1\,s^{1/\rho}}s^{\frac{1}{\rho}+1}\,\big[\frac{1}{\rho^2}s^{1/\rho}\ln
s+\ln t_0\big] ds,
$$
or
$$
|g_{2+}'(\rho)|\leq\frac{C}{|\lambda|^3
t_0^{2\rho+1}}\,\big[\frac{1}{\rho}+\ln t_0\big],
$$
where the constant $C$ depends only on $\xi^0$.

The function $g'_{2-}(\rho)$ has exactly the same estimate.

Now consider the function $g_{21}(\rho)$. For its derivative we
have
\[
g'_{21}(\rho)=\frac{a}{2\pi i \lambda^2
t_0^{\rho+1}}\cdot\int\limits_{-1}^1\frac{e^{e^{ias}}\,
e^{ias}\,e^{2ia\rho s}\big[2ias-\ln t_0-\frac{ias e^{ia\rho
s}+\lambda t_0^\rho \ln t_0}{e^{ia\rho s}+\lambda
t_0^\rho}\big]}{e^{ia\rho s}+\lambda t_0^\rho} ds.
\]
Therefore,
\[
|g'_{21}(\rho)|\leq C\,\frac{\ln t_0}{|\lambda|^3 t_0^{2\rho+1}}.
\]
Taking into account estimate (\ref{f11}) and the estimates for
$g'_{2\pm}$ and $g'_{21}$, we have
\[
sign(\lambda_1^2-\lambda_2^2)\Re(e'_{2,\lambda}(\rho))\leq
-\frac{|\lambda_1^2-\lambda_2^2|}{|\lambda|^4
t_0^{\rho+1}}+C\frac{1/\rho+\ln t_0}{|\lambda|^3 t_0^{2\rho+1}}.
\]
In other words, the left hand side is negative if
\[
t_0^{\beta_0}>\frac{C|\lambda|}{|\lambda_1^2-\lambda_2^2|}\big(\frac{1}{\beta_0}+
\ln t_0\big).
\]
Hence, there exists a number $T_1= T_1 (\xi^0, \beta_0)>1$ (see
\eqref{t0}), such, that for all $t_0\geq T_1$ the function
$sign(\lambda_1^2-\lambda_2^2)\Re(e'_{2,\lambda}(\rho))$ is
negative.

Lemma \ref{ERL} and therefore Theorem \ref{thm_Riemann-Liouville} are
proved.

\section{An example}
To illustrate the theorems proved above consider  the following
Cauchy problem (see \cite{UAY})

\begin{align} \label{ex_01}
D_{\ast}^{\beta_1} u_1(t,x) & = -D^2  u_1(t,x)-D  u_2(t,x), \quad
t>0, \ -\infty <x<\infty,
\\\label{ex_02}
D_{\ast}^{\beta_2} u_2(t,x) & = -D u_1(t,x)-D^2  u_2(t,x), \quad
t>0, \ -\infty<x<\infty,
\\ \label{ex_03}
u_1(0,x) & =\varphi_1(x), \quad u_2(0,x)=\varphi_2(x),  \quad
-\infty<x<\infty.
\end{align}
It is not hard to see that the symbol of the operator on the right
hand side of \eqref{ex_01}-\eqref{ex_02} is symmetric and has the
representation

\begin{equation} \label{matrix_10}
\mathcal{A}(\xi) = \begin{bmatrix}
        -\xi^2 & -\xi \\
        -\xi &  -\xi^2
        \end{bmatrix}
        =
        \begin{bmatrix}
        1/2 & 1/2 \\
        -1/2 &  1/2
        \end{bmatrix}
        \begin{bmatrix}
        -\xi^2+\xi & 0 \\
        0 &  -\xi^2-\xi
        \end{bmatrix}
            \begin{bmatrix}
        1 & -1 \\
        1 &  1
        \end{bmatrix}   .
        \end{equation}
As is seen from \eqref{matrix_10} that $\lambda_1(\xi)=-\xi^2+\xi$
and $\lambda_2(\xi)=-\xi^2-\xi.$ The solution   $U(t,x)= \langle
u_1(t,x), u_2(t,x) \rangle$ to Cauchy problem
\eqref{ex_01}-\eqref{ex_03} has the representation:

\begin{align*}
u_{1} (t,x)&= \frac{1}{2\pi} \int\limits_{-\infty}^{\infty}
\left[\frac{1}{2}
E_{\beta_1}((-\xi^2+\xi)t^{\beta_1})+\frac{1}{2}E_{\beta_2}((-\xi^2-\xi)t^{\beta_2})
\right] \hat{\varphi}_1](\xi) d \xi
\\
&+ \frac{1}{2\pi} \int\limits_{-\infty}^{\infty} \left[\frac{1}{2}
E_{\beta_1}((-\xi^2+\xi)t^{\beta_1})-\frac{1}{2}E_{\beta_2}((-\xi^2-\xi)t^{\beta_2})
\right] \hat{\varphi}_2](\xi) d\xi;
\\
u_{2} (t,x)&= \frac{1}{2\pi} \int\limits_{-\infty}^{\infty}
\left[\frac{1}{2}
E_{\beta_1}((-\xi^2+\xi)t^{\beta_1})-\frac{1}{2}E_{\beta_2}((-\xi^2-\xi)t^{\beta_2})
\right] \hat{\varphi}_1(\xi) d \xi
\\
&+ \frac{1}{2 \pi} \int\limits_{-\infty}^{\infty}
\left[\frac{1}{2}
E_{\beta_1}((-\xi^2+\xi)t^{\beta_1})+\frac{1}{2}E_{\beta_2}((-\xi^2-\xi)t^{\beta_2})
\right] \hat{\varphi}_2(\xi) d \xi.
        \end{align*}
Moreover, obviously, $\lambda_k(\xi) \le 0, \ k=1,2,$ for all
$\xi$ satisfying the inequality $|\xi| \ge 1.$
It is not hard to verify, that
\[
K_{1,1}(\xi, \hat{\Phi}(\xi))=K_{2,1}(\xi,
\hat{\Phi}(\xi))=\frac{1}{2}\hat{\varphi}_1(\xi)+\frac{1}{2}\hat{\varphi}_2(\xi)
\]
and
\[
K_{1,2}(\xi,
\hat{\Phi}(\xi))=\frac{1}{2}\hat{\varphi}_1(\xi)-\frac{1}{2}\hat{\varphi}_2(\xi),
\quad K_{2,2}(\xi,
\hat{\Phi}(\xi))=-\frac{1}{2}\hat{\varphi}_1(\xi)+\frac{1}{2}\hat{\varphi}_2(\xi).
\]
Therefore for the corresponding determinant one has
\[
\big|K_{j,l}\big(\xi^0,
\widehat{\Phi}(\xi^0)\big)\big|=\frac{1}{2}
\big(\hat{\varphi}^2_2(\xi^0)-\hat{\varphi}^2_1(\xi^0)\big)
\]
and condition (\ref{xi0}) has the form
\[
\hat{\varphi}^2_2(\xi^0)\neq\hat{\varphi}^2_1(\xi^0), \quad
\text{or} \quad |\hat{\varphi}_2(\xi^0)|\neq
|\hat{\varphi}_1(\xi^0)|.
\]
In this case the unknown orders $\beta_1$ and $\beta_2$ are the unique roots of the following equations
\[
E_{\beta_1} \big(\lambda_1(\xi^0) t^{\beta_1}_0\big)=
-\frac{d_1+d_2}{\hat{\varphi}_1(\xi^0)+\hat{\varphi}_2(\xi^0)},
\]
\[
E_{\beta_2} \big(\lambda_2(\xi^0) t^{\beta_2}_0\big)=
\frac{d_1-d_2}{\hat{\varphi}_2(\xi^0)-\hat{\varphi}_1(\xi^0)},
\]
respectively.

\bibliographystyle{amsplain}

\end{document}